\documentclass[12pt]{article}
\usepackage{amsmath,amsthm,amssymb,amscd,}
\usepackage{latexsym}
\newtheorem{thm}{Theorem}[section]
 \newtheorem{cor}[thm]{Corollary}

 \theoremstyle{definition}
 
 \theoremstyle{remark}

 \numberwithin{equation}{section}

\title
{A matrix differential Harnack estimate for a class of ultraparabolic equations
\thanks{Research supported by NSFC
(No. 11171025).} }

\author{ Hong Huang
  }

\date{}
\begin{document}
\maketitle

\begin{abstract}  Let $u$ be a positive solution of the ultraparabolic equation
\begin{equation*}
\partial_t u=\sum_{i=1}^n \partial_{x_i}^2 u+\sum_{i=1}^k x_i\partial_{x_{n+i}}u  \hspace{8mm}  \mbox{on} \hspace{4mm} \mathbb{R}^{n+k}\times (0,T),
\end{equation*}
where $1\leq k\leq n$ and $0<T \leq +\infty$.  Assume that $u$ and its derivatives (w.r.t. the space variables) up to the second order are bounded on any compact subinterval of $(0,T)$.  Then  the difference  $H(\log u)- H(\log f)$ of the Hessian matrices of $\log u$ and of $\log f$ (both w.r.t. the space variables) is non-negatively definite, where $f$ is the fundamental solution of the above equation with pole at the origin $(0,0)$.  The estimate in the case $n=k=1$ is due to Hamilton. As a corollary we get that $\Delta l+\frac{n+3k}{2t}+\frac{6k}{t^3}\geq 0$, where $l=\log u$,  and $\Delta=\sum_{i=1}^{n+k} \partial_{x_i}^2 $.

\vspace{0.1in} {\bf Key words}: ultraparabolic equation, matrix differential Harnack estimate,
maximum principle.

\vspace{0.1in} {\bf Mathematics Subject Classification 2010}: 35K70

\end{abstract}
\maketitle

\section{Introduction}
Since the seminal work of Li and Yau  \cite{ly}, there is extensive research  on differential Harnack inequalities for parabolic equations, see for example \cite{n} for a survey. In particular, Hamilton \cite{h1} obtained a remarkable  matrix differential Harnack estimate for the Ricci flow, whose trace form is very important. Note that  the trace version  of Hamilton's Harnack estimate for Ricci flow is derived as a corollary of his matrix Harnack estimate, and so far there is no direct proof for it (without using the matrix estimate). (See also Cao \cite{c} for a  matrix Harnack estimate for the K$\ddot{a}$hler- Ricci flow.)  In \cite{h2} Hamilton gave a matrix Harnack estimate for the heat equation on certain Riemannian manifolds, whose trace form recovers an estimate in \cite{ly}. (See also \cite{cn} for a related estimate. Recently Ni and his cooperator have further work in this direction.) This matrix Harnack estimate for the heat equation is useful for deriving monotonicity formulas, see for example \cite{h3}.

Recently, in  \cite{h4}, among other things, Hamilton extended his matrix Harnack estimate in \cite{h2} to the simple ultraparabolic equation $f_t+xf_y=f_{xx}$.
In this note we'll generalize this estimate of Hamilton in \cite{h4} to the following slightly more general class of ultraparabolic equations

\begin{equation}
\partial_t u=\sum_{i=1}^n \partial_{x_i}^2 u+\sum_{i=1}^k x_i\partial_{x_{n+i}}u  \hspace{8mm}  \mbox{on} \hspace{4mm} \mathbb{R}^{n+k}\times (0,T),
\end{equation}
where $1\leq k\leq n$ and $0<T \leq +\infty$.

(1.1) is among a still more general class of ultraparabolic equations of Kolmogorov type  satisfying the H$\ddot{o}$rmander condition (\cite{ho}); for some of the work on these equations, see for example \cite{lpp} and the references therein.  Harnack inequalities for positive solutions of these and some similar equations are extensively  studied by Polidoro  et al., see for example \cite{lp}, \cite{pp},  and  more recently, \cite{cpp} and \cite{cnp}.

The main motivation of our research is to find more  matrix differential estimates,  whose power is partially indicated above, and to  pursue more similar properties that  the heat equation shares with the  Kolmogorov type equations satisfying the H$\ddot{o}$rmander condition, which are partially displayed in some of the references cited above and
in some papers not cited here.

Using  H$\ddot{o}$rmander \cite{ho},  and Lanconelli- Polidoro \cite{lp}, one finds that the fundamental solution  of the equation (1.1) with pole at the origin (0,0) is

\vspace{0.3cm}

$\begin{array}{l}
f(x,t)=\frac{C}{t^{\frac{n+3k}{2}}}e^{-\frac{1}{t}(\sum_{i=1}^kx_i^2+\frac{1}{4}\sum_{i=k+1}^nx_i^2)
-\frac{3}{t^2}\sum_{i=1}^kx_ix_{n+i}-\frac{3}{t^3}\sum_{i=1}^kx_{n+i}^2}\\
\ \ \ \ \ \ \ \ \ =\frac{C}{t^{\frac{n+3k}{2}}}e^{-\frac{1}{4t}\sum_{i=1}^nx_i^2
-\frac{3}{t^3}\sum_{i=1}^k(x_{n+i}+\frac{1}{2}tx_i)^2},
\end{array}$

\vspace{0.3cm}

\noindent where $C$ is a constant depending only on $n$ and $k$.  Then using \cite{ho}, \cite{lp} again (see the formula (1.6) in \cite{lp}) one can easily derive the fundamental solution
\begin{equation*}
\Gamma (x,t; \xi,\tau)=\frac{C}{(t-\tau)^{\frac{n+3k}{2}}}e^{-\frac{1}{4(t-\tau)}\sum_{i=1}^n(x_i-\xi_i)^2
-\frac{3}{(t-\tau)^3}\sum_{i=1}^k(x_{n+i}-\xi_{n+i}+\frac{1}{2}(x_i+\xi_i)(t-\tau))^2}
  \end{equation*}
  with pole at any point $(\xi,\tau)$ from $f$, where $t>\tau$;  we let $\Gamma (x,t; \xi,\tau)=0$ when $t\leq \tau$.

Now let $l=$ log $f$  and   $l_{x_ix_j}=\partial_{x_i}\partial_{x_j}l$. Then

\vspace{0.3cm}

$\begin{array}{l}
  l_{x_ix_i}=-\frac{2}{t}, \hspace{2mm}  l_{x_ix_{n+i}}=l_{x_{n+i}x_i}=-\frac{3}{t^2},  \hspace{2mm} l_{x_{n+i}x_{n+i}}=-\frac{6}{t^3}  \hspace{4mm} \mbox{for}  \hspace{2mm} 1\leq i \leq k, \\
  l_{x_ix_i}=-\frac{1}{2t}   \hspace{4mm}  \mbox{for}  \hspace{2mm} k+1\leq i \leq n,  \hspace{4mm} \mbox{and}   \\
  l_{x_ix_j}=0    \hspace{4mm}  \mbox{for} \hspace{2mm} \mbox{all} \hspace{2mm} \mbox{other} \hspace{2mm} i,j.\\
\end{array}$

\vspace{0.3cm}
\noindent  Thus we get the Hessian matrix $H(\log f)=(l_{x_ix_j})_{i,j=1,\cdot\cdot\cdot,n+k}$ of log $f$ w.r.t. the space variables. Note that the matrix $((\log \Gamma (x,t; \xi,0))_{x_ix_j})=((\log f(x,t))_{x_ix_j})$ for any $\xi \in \mathbb{R}^{n+k}$.

Then  we consider a general positive solution $u$ of the equation (1.1) with $1\leq k\leq n$,
and the Hessian matrix of log $u$ w.r.t the space variables:  $H(\log u)=((\log u)_{x_ix_j})_{i,j=1,\cdot\cdot\cdot,n+k}$.

\vspace{0.3cm}

We'll use the maximum principle to show the following

\begin{thm} \label{thm 1.1} Let $u$ be a positive solution to the equation (1.1) with $1\leq k\leq n$ and $0<T \leq +\infty$. Assume that $u$ and its derivatives (w.r.t. the space variables) up to the second order are bounded on any compact subinterval of $(0,T)$. Then the Hessian $H(\log u)\geq H(\log f)$, that is, the matrix $H(\log u)- H(\log f)$ is non-negatively definite. Here, $f$ is the fundamental solution at the origin as above.
\end{thm}

\noindent This extends Theorem 4.1 in \cite{h4} which treats the case  $n=k=1$. (We also slightly weaken the assumption of Theorem 4.1 in \cite{h4}, where the solution is assumed to be bounded with bounded derivatives for $t\geq 0$. Of course  the fact that the assumption can be weakened in this way should be known to Hamilton, although he did not state it explicitly there.) The estimate in Theorem 1.1 is sharp since the equality holds trivially when $u=f$. (Note that the assumption in Theorem 1.1 on $u$ is satisfied by the fundamental solution $f$.)
This  matrix estimate contains much information. It implies that all the  principal submatrices  of the matrix $H(\log u)- H(\log f)$ are non-negatively definite. In particular, we can get some control of $l_{x_{n+i}x_{n+i}}$ in the form $l_{x_{n+i}x_{n+i}}+\frac{6}{t^3}\geq 0$ for $i=1,\cdot\cdot\cdot,k$, where $l=\log u$, even though these second order derivatives do not appear in the equation (1.1).
Below we will give three more consequences, two of which were known before (see \cite{pp}), one of which may be new.  First by partially tracing the above estimate we recover a special case of Proposition 4.2 in Pascucci- Polidoro \cite{pp}.

\begin{cor} \label{cor 1.2} (\cite{pp}) With the same assumption as in Theorem 1.1 and letting $l=$log $u$, we have

\begin{equation*}
\sum_{i=1}^nl_{x_ix_i}+\frac{n+3k}{2t}\geq 0.
\end{equation*}
\end{cor}

The original proof of Proposition 4.2 in \cite{pp}  uses a representation formula for positive solutions  of a class of Kolmogorov equations more general than (1.1). By integrating the  estimate in Corollary 1.2 along some optimal path we recover a special case of Corollary 1.2 in  \cite{pp} (see Theorem 1.2 in \cite{cpp} for an even more general version).

\begin{cor} \label{cor 1.3} (\cite{pp}) With the same assumption as in Theorem 1.1,
for any points $(p_1,\cdot\cdot\cdot,p_{n+k},t_1)$ and $(q_1,\cdot\cdot\cdot,q_{n+k},t_2)$ with $0<t_1<t_2<T$  there holds
\vspace{0.3cm}

$\begin{array}{l}
 u(q_1,\cdot\cdot\cdot,q_{n+k},t_2)\geq  \\
 (\frac{t_1}{t_2})^{\frac{n+3k}{2}}e^{-\sum_{i=1}^n\frac{(q_i-p_i)^2}{4(t_2-t_1)}-\frac{3}{(t_2-t_1)^3}\sum_{i=1}^k[q_{n+i}-p_{n+i}+\frac{1}{2}(q_i+p_i)(t_2-t_1)]^2}u(p_1,\cdot\cdot\cdot,p_{n+k},t_1).
\end{array}$

\end{cor}
\noindent Comparing the fundamental solution $\Gamma(x,t;\xi,0)$  above, one sees that the estimate in Corollary 1.3  is sharp. (This was already observed in \cite{pp}.)
This estimate in the  case  $n=k=1$ also appeared in \cite{h4}, see Corollary 4.2 there. (By the way, note that there are some misprints in the statement of Corollary 4.2 and some other places  in \cite{h4}.)
By fully tracing the matrix estimate in Theorem 1.1 we get

  \begin{cor} \label{cor 1.4}  With the same assumption as in Theorem 1.1 and letting  $l=$log $u$, there holds

\begin{equation*}
\Delta l+\frac{n+3k}{2t}+\frac{6k}{t^3}\geq 0,
\end{equation*}
where $\Delta=\sum_{i=1}^{n+k} \partial_{x_i}^2 $.
\end{cor}

  \noindent  This corollary seems to be new. It is also  sharp. Compare a similar estimate in \cite{ly} for the heat equation (see Theorem 1.1 there). Note that, as already said above, the second order derivatives $u_{x_{n+i}x_{n+i}}$ (and $l_{x_{n+i}x_{n+i}}$), for $i=1, \cdot \cdot \cdot, k$, do not appear  in the equation (1.1). For this reason,  it may not be easy to recover Corollary 1.4 by using the method in \cite{pp}. It may also be difficult to derive Corollary 1.4 by applying the maximum principle  to the scalar equation satisfied by $\Delta l$, instead of the matrix equation satisfied by the Hessian $H(l)$, since the  scalar equation satisfied by $\Delta l$ contains terms involving $l_{x_ix_j}$ for some $i\neq j$, as can be seen by tracing the equation (2.2) in  Section 2 below.  This may be another evidence for the advantage of the matrix estimates.

While the equation (1.1) is very special, we expect  similar matrix differential Harnack estimates should hold for a more general class of ultraparabolic equations of  Kolmogorov type satisfying the H$\ddot{o}$rmander condition. See Section 4 for a more precise statement.

           In the next two sections we'll prove Theorem 1.1 and Corollary 1.3 respectively, following \cite{h4} with some necessary modifications.  In Section 4 we state two conjectures related to our results above.

\section{Proof of  Theorem 1.1}

    We may and will assume that $T<+\infty$. First we claim that we can  reduce the proof of Theorem 1.1 to the case that the positive solution $u$ and its derivatives (w.r.t. the space variables) up to the second order are uniformly bounded on  $\mathbb{R}^{n+k}\times (0,T)$.  The proof of this claim is an application of a standard trick:  Suppose that $u$ is a positive solution of the equation (1.1) such that $u$ and its derivatives (w.r.t. the space variables) up to the second order are bounded on any compact subinterval of $(0,T)$. Fix $(x_0,t_0)\in \mathbb{R}^{n+k}\times (0,T)$. Let $v_\varepsilon(x,t)=u(x,t+\varepsilon)$ for $0< \varepsilon <\frac{T-t_0}{2}$.  Then $v_\varepsilon$ is a positive solution of the same equation  on $\mathbb{R}^{n+k}\times (-\varepsilon,T-\varepsilon)$ such that $v_\varepsilon$ and its derivatives (w.r.t. the space variables) up to the second order are uniformly bounded on  $\mathbb{R}^{n+k} \times [0,T-2\varepsilon]$.  If we have proven Theorem 1.1 in the case that the solution and its derivatives (w.r.t. the space variables) up to the second order are uniformly bounded, then  the matrix Harnack estimate holds for $v_\varepsilon$ at the point $(x_0,t_0)$. Note that the conclusion of the matrix Harnack estimate is independent of the bounds of $v_\varepsilon$ and its derivatives (w.r.t. the space variables) up to the second order. Then letting $\varepsilon \rightarrow 0$ one sees that the matrix Harnack inequality also holds for $u$ at $(x_0,t_0)$.

    So in the proof below we assume that the positive solution $u$ and its derivatives (w.r.t. the space variables) up to the second order are uniformly bounded on  $\mathbb{R}^{n+k}\times (0,T)$. Let $l=\log u$, and $M=H(\log u)- H(\log f)=H(l)-H(\log f)$, where $f$ is the fundamental solution of (1.1) at the origin (see Section 1). We decompose the matrix $M$ into blocks:

\vspace{0.3cm}

$M=\left(
  \begin{array}{cc}
    M_1 & M_2 \\
    M_3 & M_4 \\
  \end{array}
\right), $

\vspace{0.3cm}

\noindent where  the $n\times n$ matrix

\vspace{0.3cm}

$M_1=\left(
  \begin{array}{cccccc}
    l_{x_1x_1}+\frac{2}{t} & \cdot\cdot\cdot & l_{x_1x_k} & l_{x_1x_{k+1}} & \cdot\cdot\cdot & l_{x_1x_n}  \\
    \cdot & \cdot\cdot\cdot & \cdot & \cdot & \cdot\cdot\cdot & \cdot  \\
    l_{x_kx_1} & \cdot\cdot\cdot & l_{x_kx_k}+\frac{2}{t} & l_{x_kx_{k+1}} & \cdot\cdot\cdot & l_{x_kx_n} \\
    l_{x_{k+1}x_1} & \cdot\cdot\cdot & l_{x_{k+1}x_k} & l_{x_{k+1}x_{k+1}}+\frac{1}{2t} &\cdot\cdot\cdot & l_{x_{k+1}x_n} \\
    \cdot & \cdot\cdot\cdot & \cdot & \cdot & \cdot\cdot\cdot & \cdot  \\
    l_{x_nx_1} & \cdot\cdot\cdot & l_{x_nx_k} & l_{x_nx_{k+1}} & \cdot\cdot\cdot & l_{x_nx_n}+\frac{1}{2t} \\
  \end{array}
\right),$

\noindent  the $n\times k$ matrix

\vspace{0.3cm}

$M_2=\left(
  \begin{array}{ccc}
    l_{x_1x_{n+1}}+\frac{3}{t^2} & \cdot\cdot\cdot & l_{x_1x_{n+k}} \\
    \cdot & \cdot\cdot\cdot & \cdot \\
     l_{x_kx_{n+1}} & \cdot\cdot\cdot & l_{x_kx_{n+k}}+\frac{3}{t^2} \\
    l_{x_{k+1}x_{n+1}} & \cdot\cdot\cdot & l_{x_{k+1}x_{n+k}} \\
    \cdot & \cdot\cdot\cdot & \cdot \\
     l_{x_nx_{n+1}}  & \cdot\cdot\cdot & l_{x_nx_{n+k}} \\
  \end{array}
\right)$,

\vspace{0.3cm}

\noindent    the $k\times n$ matrix

\vspace{0.3cm}

$M_3=M_2^T, $

\vspace{0.3cm}

 \noindent and the $k\times k$ matrix

\vspace{0.3cm}

$M_4=\left(
  \begin{array}{ccc}
    l_{x_{n+1}x_{n+1}}+\frac{6}{t^3} & \cdot\cdot\cdot & l_{x_{n+1}x_{n+k}} \\
    \cdot & \cdot\cdot\cdot & \cdot \\
    l_{x_{n+k}x_{n+1}} & \cdot\cdot\cdot & l_{x_{n+k}x_{n+k}}+\frac{6}{t^3} \\
  \end{array}
\right).$

\vspace{0.3cm}

 By a direct computation  one sees that $l$ satisfies the equation
\begin{equation}
l_t=\sum_{i=1}^n(l_{x_ix_i}+l_{x_i}^2)+\sum_{i=1}^kx_il_{x_{n+i}},
\end{equation}
and $M$ satisfies the  equation
\begin{equation}
M_t=\sum_{i=1}^n(M_{x_ix_i}+2l_{x_i}M_{x_i})+\sum_{i=1}^kx_iM_{x_{n+i}}+N,
\end{equation}
where $N$ is some matrix; actually $N$ is the matrix obtained from $\tilde{N}$ below (see (2.3) and below) by replacing $\tilde{l}$ by $l$ and setting $\sigma=0, \alpha=1, \beta=1, \gamma=1$.
We want to use the maximum principle to show that the matrix $M$ is non-negatively definite. But to deal with the noncompact situation we need to apply the maximum principle to a slightly modified equation (see (2.3) below), instead of the  equation (2.2) above. So  we modify the solution $u$  to

\begin{equation*}
\tilde{u}=u+\varepsilon \{t^2\sum_{i=1}^kx_i^2+\sum_{i=1}^{n+k}x_i^2+2t(\sum_{i=1}^kx_ix_{n+i}+n)+\frac{2k}{3}t^3\}
\end{equation*}

\noindent  with $\varepsilon$ a small positive constant,  which is also a positive solution.

Let
$\tilde{l}=$ log $\tilde{u}.$
 Then  $\tilde{l}_{x_ix_j}\rightarrow 0$, for $i, j=1, \cdot\cdot\cdot, n+k$, as $|x|\rightarrow \infty$  uniformly in $t$, since now we are assuming that  $u$ and its derivatives (w.r.t. the space variables) up to the second order are uniformly bounded on  $\mathbb{R}^{n+k}\times (0,T)$.
Note that $\tilde{l} $ satisfies the equation
\begin{equation*}
\tilde{l}_t=\sum_{i=1}^n(\tilde{l}_{x_ix_i}+\tilde{l}_{x_i}^2)+\sum_{i=1}^kx_i\tilde{l}_{x_{n+i}}.
\end{equation*}

Let  $\alpha=1+\sigma\delta_0, \beta=1+\sigma\theta_0, \gamma=1+\sigma\eta_0$, where $\sigma$ is a small positive constant, and $\delta_0, \theta_0, \eta_0$ are constants  which will be chosen later, and let

\vspace{0.3cm}

$\tilde{M}=\left(
  \begin{array}{cc}
    \tilde{M}_1 & \tilde{M}_2 \\
    \tilde{M}_3 & \tilde{M}_4 \\
  \end{array}
\right), $

\vspace{0.3cm}

\noindent where  the $n\times n$ matrix

\vspace{0.3cm}

$\tilde{M}_1=\left(
  \begin{array}{cccccc}
    \tilde{l}_{x_1x_1}+\frac{2\alpha}{t} & \cdot\cdot\cdot & \tilde{l}_{x_1x_k} & \tilde{l}_{x_1x_{k+1}} & \cdot\cdot\cdot & \tilde{l}_{x_1x_n}  \\
    \cdot & \cdot\cdot\cdot & \cdot & \cdot & \cdot\cdot\cdot & \cdot  \\
    \tilde{l}_{x_kx_1} & \cdot\cdot\cdot & \tilde{l}_{x_kx_k}+\frac{2\alpha}{t} & \tilde{l}_{x_kx_{k+1}} & \cdot\cdot\cdot & \tilde{l}_{x_kx_n} \\
    \tilde{l}_{x_{k+1}x_1} & \cdot\cdot\cdot & \tilde{l}_{x_{k+1}x_k} & \tilde{l}_{x_{k+1}x_{k+1}}+\frac{1+\sigma}{2t} &\cdot\cdot\cdot & \tilde{l}_{x_{k+1}x_n} \\
    \cdot & \cdot\cdot\cdot & \cdot & \cdot & \cdot\cdot\cdot & \cdot  \\
    \tilde{l}_{x_nx_1} & \cdot\cdot\cdot & \tilde{l}_{x_nx_k} & \tilde{l}_{x_nx_{k+1}} & \cdot\cdot\cdot & \tilde{l}_{x_nx_n}+\frac{1+\sigma}{2t} \\
  \end{array}
\right),$

\vspace{0.3cm}

\noindent  the $n\times k$ matrix

\vspace{0.3cm}

$\tilde{M}_2=\left(
  \begin{array}{ccc}
    \tilde{l}_{x_1x_{n+1}}+\frac{3\beta}{t^2} & \cdot\cdot\cdot & \tilde{l}_{x_1x_{n+k}} \\
    \cdot & \cdot\cdot\cdot & \cdot \\
     \tilde{l}_{x_kx_{n+1}} & \cdot\cdot\cdot & \tilde{l}_{x_kx_{n+k}}+\frac{3\beta}{t^2} \\
    \tilde{l}_{x_{k+1}x_{n+1}} & \cdot\cdot\cdot & \tilde{l}_{x_{k+1}x_{n+k}} \\
    \cdot & \cdot\cdot\cdot & \cdot \\
     \tilde{l}_{x_nx_{n+1}}  & \cdot\cdot\cdot & \tilde{l}_{x_nx_{n+k}} \\
  \end{array}
\right)$,

\vspace{0.3cm}

\noindent    the $k\times n$ matrix

\vspace{0.3cm}

$\tilde{M}_3=\tilde{M}_2^T, $

\vspace{0.3cm}

 \noindent and the $k\times k$ matrix

\vspace{0.3cm}

$\tilde{M}_4=\left(
  \begin{array}{ccc}
    \tilde{l}_{x_{n+1}x_{n+1}}+\frac{6\gamma}{t^3} & \cdot\cdot\cdot & \tilde{l}_{x_{n+1}x_{n+k}} \\
    \cdot & \cdot\cdot\cdot & \cdot \\
    \tilde{l}_{x_{n+k}x_{n+1}} & \cdot\cdot\cdot & \tilde{l}_{x_{n+k}x_{n+k}}+\frac{6\gamma}{t^3} \\
  \end{array}
\right).$

\vspace{0.3cm}

\noindent {\bf Claim 1}  \hspace{4mm} There is a positive constant $\sigma_0$ such that for any $0<\sigma< \sigma_0$,  $\tilde{M}$ is positive definite  when $t>0$ is sufficiently small.

\vspace{0.3cm}

\noindent Claim 1 follows from the fact that there is a positive constant $\sigma_0$ such that for any $0<\sigma< \sigma_0$,  all the leading principal minors of $\tilde{M}$ are positive  when $t>0$ is sufficiently small.
The fact itself can be shown by a direct check: For the $i$-th leading principal minor, where $1\leq i \leq n$, it is trivial since now $\tilde{l}_{x_jx_l}$ (for $j, l=1,\cdot\cdot\cdot,n+k$) are uniformly bounded; for the $(n+i)$-th leading principal minor, where $1 \leq i \leq k$, it follows from that for  $t>0$ sufficiently small, the term

\vspace{0.3cm}

$ \left |
\begin{array}{ccc}
\frac{2}{t}I_k   & 0                   & \frac{3}{t^2}J \\
0                & \frac{1}{2t}I_{n-k} & 0                \\
\frac{3}{t^2}J^T & 0                   & \frac{6}{t^3}I_i \\
\end{array}
\right |$

\vspace{0.2cm}

$= \left |
\begin{array}{ccc}
\mbox{diag}(\frac{1}{2t}, \cdot\cdot\cdot,\frac{1}{2t}, \frac{2}{t}, \cdot\cdot\cdot, \frac{2}{t})   & 0                   & \frac{3}{t^2}J \\
0                & \frac{1}{2t}I_{n-k} & 0                \\
0 & 0                   & \frac{6}{t^3}I_i \\
\end{array}
\right |$

\vspace{0.2cm}

$=(\frac{1}{2t})^i(\frac{2}{t})^{k-i}(\frac{1}{2t})^{n-k}(\frac{6}{t^3})^i$

\vspace{0.3cm}

\noindent dominates the other terms  in the expansion of the $(n+i)$-th leading principal minor of $\tilde{M}$, where the $k \times i$ matrix
$J=\left(
  \begin{array}{c}
    I_i \\
    0  \\
  \end{array}
\right)$, and $I_i$ is the $i \times i$ identity matrix.

Now we choose   $(\delta_0, \theta_0$,  $\eta_0)^T$ with  $2\theta_0 \geq \eta_0$ to be an eigenvector of the matrix

\vspace{0.3cm}

$C_0:=\left(
  \begin{array}{ccc}
    -8 & 10 & -3 \\
    10 & -14 & 5 \\
    -3 & 5 & -2 \\
  \end{array}
\right)$

\vspace{0.3cm}
\noindent corresponding to a positive eigenvalue. Note that the  matrix $C_0$ does have a positive eigenvalue since its determinant  is 2.

Let  $F:=(4\alpha^2-\alpha-3\beta)(\beta^2-\gamma)-(2\alpha\beta-\beta-\gamma)^2$.  Then

\vspace{0.3cm}

$\begin{array}{l}
  F  =  (-4\delta_0^2+10\delta_0\theta_0-7\theta_0^2-3\delta_0\eta_0+5\theta_0\eta_0-\eta_0^2)\sigma^2+\cdot\cdot\cdot  \\
    \ \ \  =  \frac{1}{2}(\delta_0 \hspace{2mm}  \theta_0 \hspace{2mm} \eta_0)C_0 (\delta_0  \hspace{2mm}  \theta_0 \hspace{2mm}  \eta_0)^T\sigma^2+\cdot\cdot\cdot,
\end{array}$

\vspace{0.3cm}
\noindent where we have omitted the terms of order (w.r.t. $\sigma$) greater than 2.  Clearly there is  a positive constant $\sigma_1$ such that for any $0<\sigma< \sigma_1$, $F>0$.

\vspace{0.3cm}

\noindent {\bf Claim 2}  \hspace{4mm} With the above choice of $\delta_0, \theta_0$ and $\eta_0$, and assuming  that   $0<\sigma< $ min $\{\sigma_0, \sigma_1\}$,  the matrix $\tilde{M}$ is positive definite for all $t>0$.

\vspace{0.3cm}

\noindent  Theorem  1.1 follows from Claim 2 by first letting $\sigma\rightarrow 0$, then letting $\varepsilon\rightarrow 0$.

Before proving Claim 2 we note that $\tilde{M}$ satisfies the following equation
\begin{equation}
\tilde{M}_t=\sum_{i=1}^n(\tilde{M}_{x_ix_i}+2\tilde{l}_{x_i}\tilde{M}_{x_i})+\sum_{i=1}^kx_i\tilde{M}_{x_{n+i}}+\tilde{N},
\end{equation}
where

\vspace{0.3cm}

$\tilde{N}=\left(
  \begin{array}{cc}
    \tilde{N}_1 & \tilde{N}_2 \\
    \tilde{N}_3 & \tilde{N}_4 \\
  \end{array}
\right), $

\vspace{0.3cm}

\noindent where  the $n\times n$ matrix

\vspace{0.3cm}

$\tilde{N}_1=\left(
  \begin{array}{cc}
    \tilde{P}_1 & \tilde{P}_2 \\
    \tilde{P}_3 & \tilde{P}_4 \\
  \end{array}
\right) $

\vspace{0.3cm}

\noindent with the $k\times k$ matrix

\vspace{0.3cm}

$\tilde{P}_1=\left(                                            \begin{array}{ccc}
                                                             2\sum_{i=1}^n\tilde{l}_{x_1x_i}^2+2\tilde{l}_{x_{n+1}x_1}-\frac{2\alpha}{t^2} &\cdot\cdot\cdot & 2\sum_{i=1}^n\tilde{l}_{x_1x_i}\tilde{l}_{x_kx_i}+\tilde{l}_{x_{n+1}x_k}+\tilde{l}_{x_1x_{n+k}} \\
                                                             \cdot & \cdot\cdot\cdot & \cdot \\
                                                             2\sum_{i=1}^n\tilde{l}_{x_kx_i}\tilde{l}_{x_1x_i}+\tilde{l}_{x_{n+k}x_1}+\tilde{l}_{x_kx_{n+1}}& \cdot\cdot\cdot & 2\sum_{i=1}^n\tilde{l}_{x_kx_i}^2+2\tilde{l}_{x_{n+k}x_k}-\frac{2\alpha}{t^2} \\
                                                           \end{array}
                                                         \right),$

\vspace{0.3cm}

\noindent the $k\times (n-k)$ matrix

\vspace{0.3cm}

$\tilde{P}_2=\left(                                            \begin{array}{ccc}
                                                             2\sum_{i=1}^n\tilde{l}_{x_1x_i}\tilde{l}_{x_{k+1}x_i}+\tilde{l}_{x_{n+1}x_{k+1}} &\cdot\cdot\cdot & 2\sum_{i=1}^n\tilde{l}_{x_1x_i}\tilde{l}_{x_nx_i}+  \tilde{l}_{x_{n+1}x_n} \\
                                                             \cdot & \cdot\cdot\cdot & \cdot \\
                                                             2\sum_{i=1}^n\tilde{l}_{x_kx_i}\tilde{l}_{x_{k+1}x_i}+\tilde{l}_{x_{n+k}x_{k+1}}& \cdot\cdot\cdot & 2\sum_{i=1}^n\tilde{l}_{x_kx_i}\tilde{l}_{x_nx_i}+  \tilde{l}_{x_{n+k}x_n} \\
                                                           \end{array}
                                                         \right),$

\vspace{0.3cm}

\noindent the $(n-k)\times k$ matrix

\vspace{0.3cm}

$\tilde{P}_3=\tilde{P}_2^T$,

\vspace{0.3cm}

\noindent the $(n-k)\times (n-k)$  matrix

\vspace{0.3cm}

$\tilde{P}_4=\left(                                            \begin{array}{ccc}
                                                             2\sum_{i=1}^n\tilde{l}_{x_{k+1}x_i}^2-\frac{1+\sigma}{2t^2} &\cdot\cdot\cdot & 2\sum_{i=1}^n\tilde{l}_{x_{k+1}x_i}\tilde{l}_{x_nx_i} \\
                                                             \cdot & \cdot\cdot\cdot & \cdot \\
                                                              2\sum_{i=1}^n\tilde{l}_{x_nx_i}\tilde{l}_{x_{k+1}x_i} & \cdot\cdot\cdot & 2\sum_{i=1}^n\tilde{l}_{x_nx_i}^2-\frac{1+\sigma}{2t^2} \\
                                                           \end{array}
                                                         \right),$

\vspace{0.3cm}

\noindent  the $n \times k$ matrix

\vspace{0.3cm}

$\tilde{N}_2=\left(
  \begin{array}{ccc}
    2\sum_{i=1}^n\tilde{l}_{x_1x_i}\tilde{l}_{x_{n+1}x_i}+\tilde{l}_{x_{n+1}x_{n+1}}-\frac{6\beta}{t^3} & \cdot\cdot\cdot & 2\sum_{i=1}^n\tilde{l}_{x_1x_i}\tilde{l}_{x_{n+k}x_i}+\tilde{l}_{x_{n+1}x_{n+k}} \\
    \cdot & \cdot\cdot\cdot & \cdot\\
    2\sum_{i=1}^n\tilde{l}_{x_kx_i}\tilde{l}_{x_{n+1}x_i}+\tilde{l}_{x_{n+k}x_{n+1}} & \cdot\cdot\cdot & 2\sum_{i=1}^n\tilde{l}_{x_kx_i}\tilde{l}_{x_{n+k}x_i}+\tilde{l}_{x_{n+k}x_{n+k}}-\frac{6\beta}{t^3} \\
    2\sum_{i=1}^n\tilde{l}_{x_{k+1}x_i}\tilde{l}_{x_{n+1}x_i} & \cdot\cdot\cdot & 2\sum_{i=1}^n\tilde{l}_{x_{k+1}x_i}\tilde{l}_{x_{n+k}x_i}  \\
    \cdot & \cdot\cdot\cdot & \cdot\\
    2\sum_{i=1}^n\tilde{l}_{x_nx_i}\tilde{l}_{x_{n+1}x_i} & \cdot\cdot\cdot & 2\sum_{i=1}^n\tilde{l}_{x_nx_i}\tilde{l}_{x_{n+k}x_i} \\

  \end{array}
\right),$

\vspace{0.3cm}

\noindent  the $k \times n$ matrix

\vspace{0.3cm}

$\tilde{N}_3=\tilde{N}_2^T,$

\vspace{0.3cm}

\noindent  and finally, the $k \times k$ matrix

\vspace{0.3cm}

$\tilde{N}_4=\left(                                            \begin{array}{ccc}
                                                             2\sum_{i=1}^n\tilde{l}_{x_{n+1}x_i}^2-\frac{18\gamma}{t^4} &\cdot\cdot\cdot & 2\sum_{i=1}^n\tilde{l}_{x_{n+1}x_i}\tilde{l}_{x_{n+k}x_i} \\
                                                             \cdot & \cdot\cdot\cdot & \cdot \\
                                                              2\sum_{i=1}^n\tilde{l}_{x_{n+k}x_i}\tilde{l}_{x_{n+1}x_i} & \cdot\cdot\cdot & 2\sum_{i=1}^n\tilde{l}_{x_{n+k}x_i}^2-\frac{18\gamma}{t^4} \\
                                                           \end{array}
                                                         \right).$

\vspace{0.3cm}

Now we  prove  Claim 2 by contradiction. Note that $\tilde{M}\geq cI_{n+k}$  as $|x|\rightarrow \infty$, where $c$ is a positive constant which is uniform in $t$. Suppose  for some $0<\sigma< $ min $\{\sigma_0, \sigma_1\}$ Claim 2 is not true. Fix one such $\sigma$. Then by Claim 1 and the behavior of $\tilde{M}$ as $|x|\rightarrow \infty$, there would be a smallest $t_0 \in (0,T)$ such that there exist a point  $x_0\in \mathbb{R}^{n+k}$  and
 a nonzero vector $V =(v_1,\cdot\cdot\cdot, v_{n+k})^T \in \mathbb{R}^{n+k}$ with $\tilde{M}(x_0,t_0)V=0$.  Then at the space-time point $(x_0, t_0)$,

\vspace{0.3cm}

$\begin{array}{l}
  \sum_{j=1}^{n+k}\tilde{l}_{x_ix_j}v_j=-\frac{2\alpha}{t}v_i-\frac{3\beta}{t^2}v_{n+i}  \hspace{4mm} \mbox{for}  \hspace{2mm} 1\leq i \leq k, \\
 \sum_{j=1}^{n+k}\tilde{l}_{x_ix_j}v_j=-\frac{1+\sigma}{2t}v_i  \hspace{4mm} \mbox{for} \hspace{2mm} k+1\leq i \leq n, \hspace{2mm} \mbox{and}  \\
  \sum_{j=1}^{n+k}\tilde{l}_{x_{n+i}x_j}v_j=-\frac{3\beta}{t^2}v_i-\frac{6\gamma}{t^3}v_{n+i}  \hspace{4mm} \mbox{for} \hspace{2mm} 1\leq i \leq k.
\end{array}$

\vspace{0.3cm}

It follows that

\vspace{0.3cm}

$ \begin{array}{l}
\tilde{N}(x_0,t_0)(V,V)\\
   = 2\sum_{i=1}^n(\sum_{j=1}^{n+k}\tilde{l}_{x_ix_j}v_j)^2+2 \sum_{i=1}^kv_i(\sum_{j=1}^{n+k}\tilde{l}_{x_{n+i}x_j}v_j)\\
  -\frac{2\alpha}{t^2}\sum_{i=1}^kv_i^2-\frac{12\beta}{t^3}\sum_{i=1}^kv_iv_{n+i}-\frac{1+\sigma}{2t^2}\sum_{i=k+1}^nv_i^2-\frac{18\gamma}{t^4} \sum_{i=1}^kv_{n+i}^2\\
   =2\sum_{i=1}^k(\frac{2\alpha}{t}v_i+\frac{3\beta}{t^2}v_{n+i} )^2+2\sum_{i={k+1}}^n(\frac{1+\sigma}{2t}v_i)^2-2\sum_{i=1}^kv_i(\frac{3\beta}{t^2}v_i+\frac{6\gamma}{t^3}v_{n+i}) \\
  -\frac{2\alpha}{t^2}\sum_{i=1}^kv_i^2-\frac{12\beta}{t^3}\sum_{i=1}^kv_iv_{n+i}-\frac{1+\sigma}{2t^2}\sum_{i=k+1}^nv_i^2-\frac{18\gamma}{t^4} \sum_{i=1}^kv_{n+i}^2 \\
   =2\{\frac{4\alpha^2-\alpha-3\beta}{t^2} \sum_{i=1}^kv_i^2+\frac{6(2\alpha\beta-\beta-\gamma)}{t^3}\sum_{i=1}^kv_iv_{n+i}+\frac{9(\beta^2-\gamma)}{t^4}\sum_{i=1}^kv_{n+i}^2\}\\
   +\frac{\sigma^2+\sigma}{2t^2}\sum_{i=k+1}^nv_i^2.\\
  \end{array}$

\vspace{0.3cm}

By our choice of $\delta_0, \theta_0$,  $\eta_0$  and $\sigma$,  $\beta^2-\gamma=\sigma^2\theta_0^2+\sigma(2\theta_0-\eta_0)>0$  (noting that  $2\theta_0-\eta_0$ and $\theta_0$ can not both be zero since $(1,0,0)^T$ is not an eigenvector of the matrix $C_0$), and  $F=(4\alpha^2-\alpha-3\beta)(\beta^2-\gamma)-(2\alpha\beta-\beta-\gamma)^2>0$. On the other hand, by Cauchy-Schwarz inequality $(\sum_{i=1}^kv_iv_{n+i})^2\leq \sum_{i=1}^kv_i^2\sum_{i=1}^kv_{n+i}^2$. It follows that  $\tilde{N}(x_0, t_0)(V,V)>0$.

Now we arrive at a  contradiction by applying  the equation (2.3) at $(x_0,t_0)$  to $(V,V)$: $\tilde{M}_t(x_0,t_0)(V,V)\leq 0$, $\tilde{M}_{x_i}(x_0,t_0)(V,V)=0$, $\tilde{M}_{x_ix_i}(x_0,t_0)(V,V)\geq 0$, for $1\leq i\leq n+k$, but $\tilde{N}(x_0, t_0)(V,V)>0$.  This completes the proof of Claim 2 and Theorem 1.1.

\vspace{0.3cm}

\noindent {\bf Remark}  \hspace{4mm} From the proof above we see that if we assume  that  $l_{x_ix_j}\rightarrow 0$ ($i,j=1,\cdot\cdot\cdot,n+k$, where $l=\log u$) as $|x|\rightarrow \infty$  uniformly for  $t$ in any compact subinterval of $(0,T)$, instead of assuming that $u$ and its derivatives (w.r.t. the space variables) up to the second order are bounded on any compact subinterval of $(0,T)$ as in the original statement of Theorem 1.1, then the result of Theorem 1.1 also holds true.

\section{Proof of the Corollary 1.3}

   We follow closely Hamilton \cite{h4}. From the equation (2.1) satisfied by $l$ and Corollary 1.2 we get that
\begin{equation*}
l_t\geq -\frac{n+3k}{2t}+\sum_{i=1}^nl_{x_i}^2+\sum_{i=1}^kx_il_{x_{n+i}}.
\end{equation*}
Along any path with $\frac{dx_{n+i}}{dt}=-x_i$  ($1\leq i\leq k$)  we compute
\begin{equation*}
\frac{dl}{dt}=l_t+\sum_{i=1}^{n+k}l_{x_i}\frac{dx_i}{dt}\geq \sum_{i=1}^n(l_{x_i}^2+l_{x_i}\frac{dx_i}{dt})-\frac{n+3k}{2t}\geq -\frac{n+3k}{2t}-\frac{1}{4}\sum_{i=1}^n(\frac{dx_i}{dt})^2.
\end{equation*}

We integrate along such path and get
\begin{equation}
l(q_1,\cdot\cdot\cdot,q_{n+k},t_2)\geq l(p_1,\cdot\cdot\cdot,p_{n+k}, t_1)-\frac{n+3k}{2}\log \frac{t_2}{t_1}-\frac{1}{4}\int_{t_1}^{t_2}\sum_{i=1}^n(\frac{dx_i}{dt})^2dt.
\end{equation}

The optimal path will  minimize the integral

\begin{equation*}
\int_{t_1}^{t_2}\sum_{i=1}^n(\frac{dx_i}{dt})^2dt
\end{equation*}
with the constraints that
\begin{equation*}
\int_{t_1}^{t_2}x_idt=-(q_{n+i}-p_{n+i})  \hspace{4mm} (1\leq i \leq k).
\end{equation*}
The Euler-Lagrange equations give that along the optimal path $\frac{d^2x_i}{dt^2}$ are constants independent of  $t$ for $1\leq i \leq k$, and $\frac{dx_i}{dt}$ are constants independent of  $t$ for $k+1\leq i \leq n$. So such path should  have the form

\vspace{0.3cm}

$\begin{array}{l}
  x_i=3a_it^2+2b_it+c_i,  \hspace{4mm} 1\leq i \leq k, \\
  x_i=d_it+e_i,  \hspace{4mm} k+1\leq i \leq n,\\
  x_{n+i}=-(a_it^3+b_it^2+c_it+f_i), \hspace{4mm} 1\leq i \leq k,
\end{array}$

\vspace{0.3cm}
\noindent where $a_i$, $b_i, \cdot\cdot\cdot$, $f_i$ are constants.
As in Hamilton \cite{h4} we compute the optimal path from $(p_1, \cdot\cdot\cdot, p_{n+k}, t_1)$ to $(q_1, \cdot\cdot\cdot, q_{n+k}, t_2)$, using the substitution
\begin{equation*}
x_i=\hat{x}_i+\frac{q_i-p_i}{t_2-t_1}t+\frac{p_it_2-q_it_1}{t_2-t_1}  \hspace{4mm} (1\leq i \leq n).
\end{equation*}

Now the problem is reduced to minimize
\begin{equation*}
\int_{t_1}^{t_2}\sum_{i=1}^n(\frac{d\hat{x}_i}{dt})^2dt
\end{equation*}
with the constraints that
\begin{equation*}
\int_{t_1}^{t_2}\hat{x}_idt=-(q_{n+i}-p_{n+i})-\frac{1}{2}(q_i+p_i)(t_2-t_1)  \hspace{4mm} (1\leq i \leq k),
\end{equation*}
and the boundary conditions
\begin{equation*}
\hat{x}_i=0   \hspace{4mm} \mbox{at} \hspace{4mm} t=t_1  \hspace{4mm} \mbox{and} \hspace{4mm} \mbox{at} \hspace{4mm}  t=t_2 \hspace{4mm} \mbox{for}   \hspace{4mm} 1\leq i \leq n.
\end{equation*}
The solution is given by

\vspace{0.3cm}

$\begin{array}{l}
  \hat{x}_i=\frac{6}{(t_2-t_1)^3} [-(q_{n+i}-p_{n+i})-\frac{1}{2}(q_i+p_i)(t_2-t_1)](t_2-t)(t-t_1)  \hspace{2mm} \mbox{for}   \hspace{2mm} 1\leq i \leq k,\\
  \hat{x}_i=0   \hspace{4mm} \mbox{for}   \hspace{4mm} k+1\leq i \leq n.
\end{array}$

\vspace{0.3cm}

Now
\begin{equation*}
\int_{t_1}^{t_2}\sum_{i=1}^n(\frac{d\hat{x}_i}{dt})^2dt=\frac{12}{(t_2-t_1)^3}\sum_{i=1}^k[q_{n+i}-p_{n+i}+\frac{1}{2}(q_i+p_i)(t_2-t_1)]^2,
\end{equation*}
and

\begin{equation*}
\int_{t_1}^{t_2}\sum_{i=1}^n(\frac{dx_i}{dt})^2dt=\sum_{i=1}^n\frac{(q_i-p_i)^2}{t_2-t_1}+\frac{12}{(t_2-t_1)^3}\sum_{i=1}^k[q_{n+i}-p_{n+i}+\frac{1}{2}(q_i+p_i)(t_2-t_1)]^2.
\end{equation*}

Then we insert this in (3.1) and Corollary 1.3 follows by exponentiating.

\section{Two conjectures}

First we propose

\vspace{0.3cm}

\noindent {\bf Conjecture 1 \hspace{2mm}}  Theorem 1.1 still holds true without  assuming that  $u$ and its derivatives (w.r.t. the space variables) up to the second order are bounded on any compact subinterval of $(0,T)$.

\vspace{0.3cm}

Compare the last remark in Section 2. One  (non-direct) evidence for this conjecture is that there is no such assumption in Corollary 1.2 in  \cite{pp}.  Perhaps one way to prove Conjecture 1 is to try to localize the estimate in Section 2 above.  Note that there is a localized (non-matrix) differential Harnack estimate in  Li-Yau [12] for the heat equation.  But so far, even for the heat equation, the localized matrix differential Harnack estimate is missing.

To state our second conjecture  let
\vspace{0.3cm}

$A=\left(
  \begin{array}{cc}
    A_0 & 0 \\
    0 & 0 \\
  \end{array}
\right), $

\vspace{0.3cm}

\noindent and

\vspace{0.3cm}

$B=\left(
  \begin{array}{ccccc}
    0 & B_1  & 0 & \cdot\cdot\cdot & 0  \\
    0 & 0 & B_2 &  \cdot\cdot\cdot & 0  \\
     \cdot & \cdot & \cdot & \cdot\cdot\cdot &\cdot \\
    0 & 0 & 0 & \cdot\cdot\cdot & B_r \\
    0 & 0 & 0  & \cdot\cdot\cdot & 0 \\
      \end{array}
\right),$

\vspace{0.3cm}

\noindent be two constant real $N\times N$ matrices (for some $N$), where $A_0$ is a positive definite symmetric  $p_0 \times p_0$ matrix, and $B_i$ is a $p_{i-1}\times p_i$ matrix of rank $p_i$ for $i=1,\cdot\cdot\cdot,r$, where

\vspace{0.3cm}

$p_0\geq p_1 \geq \cdot\cdot\cdot \geq p_r \geq 1 \hspace{8mm} \mbox{and} \hspace{8mm} \sum_{i=0}^{r}p_i=N$.

\vspace{0.3cm}

Then let the operator
\begin{equation*}
L=\mbox{div}(AD)+<x,BD>-\partial_t,  \hspace{8mm} (x,t)\in \mathbb{R}^N \times (0,T),
\end{equation*}
where $D=(\partial_{x_1},\cdot\cdot\cdot,\partial_{x_N})^T$, and div is the divergence in $\mathbb{R}^N$.

Note that for the corresponding operator in our equation (1.1), $N=n+k$, $r=1$, $A_0$ is the $n\times n$ identity matrix $I_n$, and $B_1$ is the $n\times k$ matrix

\vspace{0.2cm}

$\left(
  \begin{array}{c}
    I_k  \\
    0  \\
       \end{array}
\right).$

\vspace{0.2cm}

From \cite{lp} we know that the operator $L$ satisfies  H$\ddot{o}$rmander's hypoellipticity condition (\cite {ho}). In \cite{pp} (see also \cite{cpp}) a Harnack estimate is obtained for the equation $Lu=0$.
Now we propose

\vspace{0.3cm}

\noindent {\bf Conjecture 2 \hspace{2mm}}  The same result as in  Theorem 1.1 still holds for the more general equation $Lu=0$ for $L$ defined above; of course now the fundamental solution $f$ of (1.1) should be replaced by that  of the equation $Lu=0$.

\vspace{0.3cm}

If Conjecture 2 should be  true then it would recover the Harnack estimate in \cite{pp} at least when the solution  and its derivatives (w.r.t. the space variables) up to the second order are bounded on any compact subinterval of $(0,T)$.
 I hope that when $r\leq 2$  Conjecture 2 could be attacked by an argument similar to that in this paper. This will be checked in our future study.  In general I expect the geometry of the operator $L$ will play some role.

\vspace{0.3cm}

\noindent {\bf Acknowledgements} \hspace{4mm} I would like to thank the referee for helpful comments.

\vspace*{0.4mm}
 School of Mathematical Sciences,  Key Laboratory of Mathematics and
Complex Systems,

Beijing Normal University,

Beijing 100875, P.R. China

 E-mail address: hhuang@bnu.edu.cn

\end{document}